\documentclass[11pt,twoside]{article}
\usepackage{fancyhead}
\usepackage{psfrag}
\usepackage{epsfig}
\usepackage{graphicx}
\usepackage{amsmath}
\usepackage{amsthm}
\usepackage{subfigure}
\usepackage{cite}
\usepackage[latin1]{inputenc}

\begin{document}

\date{}

\vspace*{5mm}

\noindent \textbf{\LARGE Minimal complete arcs in $PG$(2,$q$), $q\leq 32$} %
\thispagestyle{fancyplain} \setlength\partopsep {0pt} \flushbottom

\vspace*{5mm}

\noindent \textsc{Stefano Marcugini} \hfill \texttt{gino@dipmat.unipg.it}%
\newline
\textsc{Alfredo Milani}\hfill \texttt{milani@unipg.it}\newline
\textsc{Fernanda Pambianco}\hfill \texttt{fernanda@dipmat.unipg.it}\newline
{\small Dipartimento di Matematica e Informatica, Università degli Studi di
Perugia,\newline
Via Vanvitelli~1, Perugia, 06123, Italy}

\medskip

\begin{center}
\parbox{11,8cm}{\footnotesize\textbf{Abstract.}
In this paper it has been verified, by a computer-based proof, that the
smallest size of a complete arc is $14$ in $PG(2,31)$ and in $PG(2,32)$.
Some examples of  such arcs are also described.
}
\end{center}

\baselineskip=0.9\normalbaselineskip

\section{Introduction}

In the projective plane $PG(2,q)$ over the Galois field $GF(q)$ an $n-$arc
is a set of $n$ points no $3$ of which are collinear. An $n$-arc is called
complete if it is not contained in an $(n+1)$-arc of the same projective
plane. For a detailed description of the most important properties of these
geometric structures, we refer the reader to \cite{librohir}. In \cite
{surveyhir} the close relationship between the theory of complete $n-$arcs,
coding theory and mathematical statistics is presented. In particular arcs
and linear maximum distance separable codes (MDS codes) are equivalent
objects (see \cite{tamas2}, \cite{tallini}, \cite{thas}). Partly because of
this fact, in recent years, the problem of determining the spectrum of
values $n$ for which a complete arc exists has been intensively
investigated. For recent results on the sizes of complete arcs in projective
planes see \cite{davArchi}. The full classification of complete $n-$arcs is
known for $q\leq 29$, see \cite{Cool29} and the references therein. This
paper concerns the minimal complete arcs in $PG(2,q)$ for $q\leq 32$. The
minimal size of a complete $n$-arc of $PG(2,q)$ is indicated by $t(2,q)$.
General lower bounds on $t(2,q)$\ are given in the following table:

\smallskip \smallskip

\begin{center}
\begin{tabular}{c}
\begin{tabular}{|c|c|c|}
\hline
$q$ & $t(2,q)>$ & References \\ \hline\hline
$q$ & $\sqrt{2q}+1$ & \cite{hirs132} \\ \hline
$q=p^{h}$, $p$ prime, $h=1,2,3$ & $\sqrt{3q}+1/2$ & \cite{hirs8}, \cite
{hirs18}, \cite{polv} \\ \hline
\end{tabular}
\medskip \\
Lower bounds for $t(2,q)$%
\end{tabular}
\end{center}

\smallskip

The values of $t(2,q),q\leq 29$ are stated in the following table:

\smallskip

\begin{center}
\begin{tabular}{c}
\begin{tabular}{|c|c|c|c|c|}
\hline
$q$ & $t(2,q)$ &
\begin{tabular}{c}
Number of classes \\
up to $PGL(2,q)$%
\end{tabular}
&
\begin{tabular}{c}
Number of classes \\
up to $P\Gamma L(2,q)$%
\end{tabular}
& References \\ \hline\hline
$2$ & $4$ & $1$ &  & \cite{librohir} \\ \hline
$3$ & $4$ & $1$ &  & \cite{librohir} \\ \hline
$4$ & $6$ & $1$ & $1$ & \cite{librohir} \\ \hline
$5$ & $6$ & $1$ &  & \cite{librohir} \\ \hline
$7$ & $6$ & $2$ &  & \cite{librohir} \\ \hline
$8$ & $6$ & $3$ & $1$ & \cite{librohir}, \cite{report23sec} \\ \hline
$9$ & $6$ & $1$ & $1$ & \cite{librohir} \\ \hline
$11$ & $7$ & $1$ &  & \cite{6}, \cite{faina24} \\ \hline
$13$ & $8$ & $2$ &  & \cite{tesi13}, \cite{faina24} \\ \hline
$16$ & $9$ & $6$ & $2$ & \cite{faina40}, \cite{report23sec}, \cite{faina59}
\\ \hline
$17$ & $10$ & $560$ &  & \cite{faina15}, \cite{faina59} \\ \hline
$19$ & $10$ & $29$ &  & \cite{faina15}, \cite{faina59} \\ \hline
$23$ & $10$ & $1$ &  & \cite{faina15} \\ \hline
$25$ & $12$ &  & $606$ & \cite{art25} \\ \hline
$27$ & $12$ &  & $7$ & \cite{art29} \\ \hline
$29$ & $13$ & $708$ &  & \cite{Cool29}, \cite{art29} \\ \hline
\end{tabular}
\medskip  \\
Value of $t(2,q),q\leq 29$%
\end{tabular}
\end{center}

\smallskip

In this paper it is demonstrated by a computer-based proof that

\begin{center}
$t_{2}(2,31)=14$ \thinspace \thinspace \thinspace \thinspace \thinspace
\thinspace \thinspace \thinspace \thinspace \thinspace and $%
\,\,\,\,\,\,\,\,\,\,\,t_{2}(2,32)=14$.
\end{center}

This result has been obtained by an exhaustive computer search. The search
has been feasible because projective equivalence properties among arcs have
been exploited and a simple parallelization technique has been used.

We also performed a partial classification of the smallest complete arcs in $%
PG(2,31)$ and in in $PG(2,32)$ obtaining $3391$ and $9300$
non-equivalent examples respectively. Equivalence up to $PGL(3,31)$
has been considered for $PG(2,31)$, while equivalence up to $P\Gamma
L(3,32)$ has been considered for $PG(2,31)$. The aim of this search
has been to look for examples with large automorphism group, but the
maximum order of the automorphism group of the examined examples is
$6$ and almost all the examples have trivial automorphism group.

In Section 2 the computation of the values $t_{2}(2,31)$ and $t_{2}(2,32)$
is described; examples of the smallest complete arcs in $PG(2,31)$ and in in
$PG(2,32)$ are presented in Section 3.

\smallskip

\section{The determination of $t_{2}(2,31)$, $t_{2}(2,32)$}

\qquad The results presented in this paper have been obtained by an
exhaustive computer search. The exhaustive search has been feasible because
projective properties among arcs have been exploited to avoid obtaining too
many isomorphic copies of the same solution arc and to avoid searching
through parts of the search space isomorphic to previously searched portions.

Also a simple parallelization technique has been used to divide the load of
the computation in a multiprocessor computer (a Quad-Core Linux computer
with 4 processors).

The used algorithm starts constructing a tree structure containing a
representative of each class of non-equivalent arcs of size less than or
equal to a fixed threshold $h$. If the threshold $h$ were equal to the
actual size of the sought arcs, the algorithm would be orderly, that is
capable of constructing each goal configuration exactly once \cite{orderly}.
However, in the present case, the construction of the tree with the
threshold $h$ equal to the size of the sought arcs would have been too space
and time consuming. For this reason a hybrid approach has been adopted. The
tree representing the non-equivalent arcs of size less than or equal to
eight has been constructed and then every non-equivalent $8-$arc has been
extended using a backtracking algorithm trying to obtain complete arcs of
the desired size. In the backtracking phase, the information obtained during
the classification of the arcs has been further exploited to prune the
search tree. In fact the points that would have given arcs equivalent to
already obtained ones have been excluded from the backtracking steps. The
algorithm is described in detail in \cite{art25}.

Each $8-$arc can be extended in independent way. To distribute the load of
computation, a certain number of 8-arcs has been assigned to each of the 4
processors to be extended with the backtracking algorithm.

As the backtracking algorithm exploits the information obtained during the
classification phase, the extension time of the $8-$arcs is not equal. The $%
8-$arcs are extended following a certain order; when we extend an
$8-$arc we can avoid to consider some possibilities because we know
that we should obtain solutions equivalent to solutions obtained
extending $8-$arcs already considered. It means that the extension
time of the first $8-$arcs is much longer of the extension time of
the following $8-$arcs. Therefore, to balance the computational load
among the 4 processors, we divided the number of $8-$arcs to extend
according to the following proportions: 10\%, 20\%, 30\%, 40\%.

When \ studying \ the \ value\ of $t_{2}(2,31)$, during the classification,
up to $PGL(3,31)$, of the arcs of $PG(2,31)$\ of size less than or equal to
8, we have found 11 non-equivalent arcs of size five, 905 non-equivalent
arcs of size six, 66,272 non-equivalent arcs of size seven and 3,768,298
non-equivalent arcs of size eight. Each 8-arc has been extended trying to
obtain complete arcs of size less than or equal to 13. No examples have been
found, so $t_{2}(2,31)=14$.

When \ studying \ the \ value\ of $t_{2}(2,32)$,\ during the classification,
up to $P\Gamma L(3,32)$, of the arcs of $PG(2,32)$\ of size less than or
equal to 8, we have found 3 non-equivalent arcs of size five, 213
non-equivalent arcs of size six, 16,593 non-equivalent arcs of size seven
and 1,031,750 non-equivalent arcs of size eight.

Each 8-arc has been extended trying to obtain complete arcs of size less
than or equal to 13. No examples have been found, so $t_{2}(2,32)=14$.

The search for the $13-$arcs in $PG(2,31)$\ lasted about 197 days of total
CPU time, while the search for the $13-$arcs in $PG(2,32)$\ lasted about 100
days of total CPU time. The reason because the search in the bigger plane
has been quicker is that $P\Gamma L(3,32)$ is much bigger than $PGL(3,31)$,
so the consideration about isomorphism properties have reduced the search
space in $PG(2,32)$, as we can see by the reduced number of classes to
extend.

\section{Examples of the smallest complete arcs in PG(2,31) and in PG(2,32)}

After having investigated the values of $t_{2}(2,31)$\thinspace \thinspace
and$\,\,t_{2}(2,32)$, we performed a partial search for examples of $14-$%
arcs using the same algorithm.

We stopped the search for $14-$arcs in $PG(2,31)$\ after having
found 500,000 examples. We performed a partial classification of
them using MAGMA, a system for symbolic computation developed at the
University of Sydney.

We obtained $3286$ non-equivalent examples with trivial stabilizer group, $%
97 $ non-equivalent examples with stabilizer group of order two, $3$
non-equivalent examples with stabilizer group isomorphic to $Z_{4}$, $4$
non-equivalent examples with stabilizer group isomorphic to $Z_{2}\times
Z_{2}$ and one example with stabilizer group isomorphic to $S_{3}$.

We have stopped the search for $14-$arcs in $PG(2,32)$\ after having
found 20,000 examples.

After a partial classification using MAGMA, we obtained $8759$
non-equivalent examples with trivial stabilizer group and $541$
non-equivalent examples with stabilizer group of order two, one
example with stabilizer group isomorphic to $Z_{4}$ and one example
with stabilizer group isomorphic to $Z_{5}$.

The field $GF(32)$ has been constructed using the primitive polynomial $\xi
^{3}+2\xi ^{2}+1$. Let $R=\{(0,0,1),(0,1,0),(1,0,0),(1,1,1)\}$.

The $14-$arc in $PG(2,31)$\ with stabilizer group $S_{3}$ is:

\begin{center}
\thinspace $K_{1}=R\cup \{(1,3,10)$, $(1,5,11)$, $(1,9,29)$, $(1,12,19)$, $%
(1,13,6)$, $(1,14,3)$, $(1,16,9)$, $(1,20,26)(1,21,15),(1,22,16)\}.$
\end{center}

The $14-$arc in $PG(2,32)$\ with stabilizer group $Z_{4}$ is:

\begin{center}
\thinspace $K_{2}=R\cup \{(1,\xi ^{17},\xi ^{2})$, $(1,\xi ^{27},\xi ^{4})$,
$(1,\xi ^{26},\xi ^{13})$, $(1,\xi ^{2},\xi )$, $(1,\xi ,\xi ^{28})$, $%
(1,\xi ^{18},\xi ^{19})$, $(1,\xi ^{15},\xi ^{9})$, $(1,\xi ^{20},\xi
^{10}),(1,\xi ^{23},\xi ^{29}),(1,\xi ^{10},\xi ^{12})\}.$
\end{center}

The $14-$arc in $PG(2,32)$\ with stabilizer group $Z_{5}$ is:

\begin{center}
\thinspace $K_{3}=R\cup \{(1,\xi ^{24},\xi )$, $(1,\xi ^{7},\xi ^{14})$, $%
(1,\xi ,\xi ^{28})$, $(1,\xi ^{8},\xi ^{18})$, $(1,\xi ^{28},\xi ^{10})$, $%
(1,\xi ^{22},\xi ^{7})$, $(1,\xi ^{2},\xi ^{22})$, $(1,\xi ^{23},\xi
^{29}),(1,\xi ^{30},\xi ^{13}),(1,\xi ^{25},\xi ^{23})\}.$
\end{center}

\end{document}